# Active Inference for Energy Control and Planning in Smart Buildings and Communities

Seyyed Danial Nazemi [1,*], Mohsen A. Jafari [1], and Andrea Matta [2]

*Abstract*— Active Inference (AIF) is emerging as a powerful framework for decision-making under uncertainty, yet its potential in engineering applications remains largely unexplored. In this work, we propose a novel dual-layer AIF architecture that addresses both building-level and community-level energy management. By leveraging the free energy principle, each layer adapts to evolving conditions and handles partial observability without extensive sensor information and respecting data privacy. We validate the continuous AIF model against both a perfect optimization baseline and a reinforcement learning-based approach. We also test the community AIF framework under extreme pricing scenarios. The results highlight the model's robustness in handling abrupt changes. This study is the first to show how a distributed AIF works in engineering. It also highlights new opportunities for privacy-preserving and uncertainty-aware control strategies in engineering applications.

## I. INTRODUCTION

Active Inferencing (AIF) is a mathematical framework rooted in the Free Energy Principle (FEP), which states that organisms maintain themselves in states in which they expect to be over a long time. An AIF agent, therefore, aims at minimizing the surprisal of the outcomes or observations that are experienced through its interactions with the world. AIF works based on two models [1], namely, Generative Density (G-density) and Recognition Density (R-density), which, respectively, describe the world (actual process) and the perception of the AIF agent. The generative process generates observations and is subject to actions (if any) prescribed by the agent. There are two options before the agent if these observations do not match the agent's prediction: (a) update its beliefs (perception model), (b) take action to change the world. The agent decides on these options using Free energy which has retrospective and prospective components, namely, Variational Free Energy (VFE) and Expected Free Energy (EFE). Minimizing VFE updates the perception model whereas minimizing EFE determines the "optimal" policy to plan. Originally formulated in neuroscience and mathematical psychology to explain perception and action in the brain, AIF has been utilized for state estimation, control, and decision-making in robotics [2]. In cognitive science, AIF has been linked to goal-directed behavior and epistemic foraging, providing insights into human decision-making and learning processes [3]. Recent works have also incorporated AIF into vehicle navigation, where AIF manages uncertainties without requiring extensive training or computational resources, offering significant advantages over traditional methods like Model Predictive Control (MPC) and Reinforcement Learning (RL). Unlike deterministic optimization methods, which struggle with inaccuracies in forecasts, and RL models, which require significant computational resources for training, AIF employs Bayesian principles to handle uncertainty and make real-time adaptive decisions. This ability to infer latent states and incorporate incomplete or noisy data makes AIF particularly suited for complex and evolving systems with non-stationary behavior. Its reliance on generative modeling reduces the need for intrusive data collection, addressing privacy concerns that are increasingly critical in modern energy systems. These features position AIF as a compelling alternative to traditional methods, enabling secure and adaptive energy management. We are motivated by a number of factors: (a) many engineering systems are not fully observable and they are subject to underlying hidden uncertainties, (b) many real-life systems are hybrid and demonstrate both time- and event-based dynamics, (c) AIF is organic in maintaining a perception model that can best describe a real process amid the observed outcomes, (d) AIF by its nature is a non-stationary model and, hence, can describe non-stationary engineering systems better than traditional models, such as RL or MPC, and (e) except for a few sporadic works in the literature, mainly in robotic control, the application of AIF to engineering problems is yet to be seen.

In this article, we demonstrate a dual-layer AIF architecture for building energy management applications. Our reasoning is that building thermal dynamics (e.g., temperature changes) are non-stationary, influenced by many internal and external factors and uncertainties, some of which are not observable nor controllable. For instance, changes in occupancy are hard to observe (in real-time) and can significantly impact the heat and humidity balance in a building zone. Similarly, unsupervised opening of a window (usually an unobservable event) exposes a zone to external air infiltration and can immediately change the heat balance within that zone. From where a community manager stands, there may be many hidden causes and uncertainties at building levels that influence their power consumption patterns. Furthermore, for privacy reasons, buildings may refrain from disclosing internal data to the community manager. AIF allows for these to be inferred from more readily observable data, such as changes in zone temperature or changes in building power consumption.

The proposed dual layer includes continuous AIF agents in Layer 1 and a discrete AIF agent for the community manager in Layer 2. A building AIF agent tracks a preferred target set by the building manager, e.g., "optimal" cooling conditions for its occupants. This is accomplished by controlling the heating, ventilation, and air conditioning (HVAC) resources amid internal hidden dynamics and uncertainties within that building. At the community layer, the discrete AIF agent tracks a preferred target, e.g., keeping the community power consumption within the proximity of a day-ahead plan. These day-ahead plans are determined using forecasts and market pricing schemes. The target tracking is accomplished by managing the peak load, common energy storage system (ESS), and distributed renewables and participating in market transactions in response to external price signals. Linking the two layers is a privacy-preserving interaction protocol, which ensures that occupant-level data—like real-time schedules or appliance usage—remain undisclosed, and only partial observations are maintained by the community. Our formulation is scalable, and computational costs are far less than the usual data-driven models, such as RL. In particular, the formulation can easily extend to many individualized continuous AIF agents, each with large sets of latent variables and control actions for a building. The community layer AIF can accommodate many continuous AIF agents under its control, and its action and latent state spaces can be quite large. The work is pioneering in its problem formulation and the AIF approach to control and planning.

[1] Seyyed Danial Nazemi and Mohsen A. Jafari are with the Department of Industrial and Systems Engineering, Rutgers, The State University of New Jersey, Piscataway, NJ 08854, USA. Email: danial.nazemi@rutgers.edu.

[2] Andrea Matta is with the Department of Mechanical Engineering, Politecnico di Milano, 20156 Milan, Italy.

Optimization-based approaches have been widely used in building energy management. These methods often rely on explicit optimization models solved through heuristic algorithms [4], mixed-integer linear programming (MILP) [5] or mixed-integer nonlinear programming (MINLP) [6]. Forecast errors and the computational complexity of multi-objective optimization pose significant challenges, limiting their performance and scalability. Model Predictive Control (MPC) is a relatively popular framework for building energy management due to its rolling horizon optimization [7]. MPCs can also be resilient to prediction errors, but their performance heavily depends on the precision of embedded models. Reinforcement Learning (RL), Multi-agent RL, and Deep RL methods have been widely applied in energy management to optimize energy consumption, reduce costs, and enhance system reliability [8]. Besides training complexity, RL-based methods face challenges with scalability and privacy due to centralized architectures and are susceptible to security vulnerabilities, such as cyber-attacks. Time-variant and non-stationary dynamics and partial information are also major issues with RL-based methods. Model-based RL can mitigate some of these issues, but planning and model updates are still isolated. As we will demonstrate, AIF integrates planning and model updates into a single framework, ensuring plans always use the most up-to-date model of the underlying process.

## II. Problem Statement and Preliminaries

In this article, the community is composed of n=1,…,N buildings, k=1,2,…,K community shared assets, the energy market, and a transactional model that encompasses energy transactions between the community and the rest of the power grid. The problem of interest is to formulate continuous AIF agents at building levels and a discrete agent at the community level to track the respective preferred targets, usually defined on outcomes and may change dynamically. Figure 1 gives a schematic of a connected community. The problem formulation consists of the general Bayesian formulation for each AIF agent and recognition of its hidden states, actions, and observations or outcomes. For continuous agents, we will define hidden states by $x \varepsilon \chi$, actions by $a \varepsilon \mathbb{A}$ and observation by $\phi \varepsilon \Phi$ with cardinalities $|\chi|, |\mathbb{A}|$ and $|\Phi|$, respectively. For discrete agents, we define states by $s \in \mathcal{S}$, actions by $u \in \mathcal{U}$, and outcomes by $o \varepsilon O$, with cardinalities $|\mathcal{S}|, |\mathcal{U}|$ and $|O|$, respectively. For building-level agents, these variables will have an index of $n$, but we will deliberately drop the index to simplify the presentation.

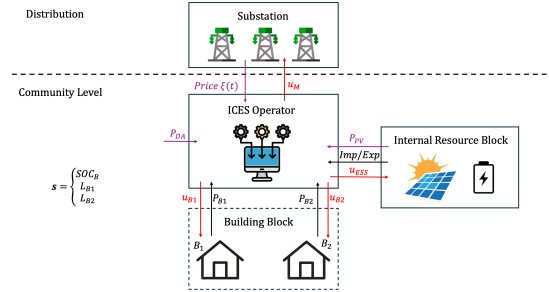

Figure 1. A schematic diagram of a connected building community with two buildings and one internal resource block connected to the power grid

For a continuous agent, we define a *likelihood function* $g(.)$, which describes how observed data is generated from hidden states and control actions, and a dynamic function $f(.)$. The latter one, if not known from the physics of the problem, is approximated or trained. For a discrete AIF agent, we must define a likelihood matrix $p(o(t)|s(t)) = A \in \mathbb{R}^{(num_{obs}) \times (num_{states})}$ and a transition matrix $p(s(t+1)|s(t), u(t) = u) = B[u(t)](s(t), s(t+1))]$ for a one-step transition depending on the action policies. The generative model of each building or the community manager defines the respective probabilistic relationship between states and outcomes; for instance, for a building agent, we have $p(\phi, x) = p(\phi|x)p(x)$. The *prior distribution*, $p(x)$, captures the initial beliefs about the dynamics of hidden states over time. Extensions of this relation can incorporate additional actions and parameters. We are particularly interested in inferring the posterior probability of hidden states given outcomes, i.e., $p(x|\phi)$. However, directly computing this distribution is often intractable, necessitating the use of **variational inference**. We will use $q(.)$ It to represent variational (or R-Density), and the Kullback-Leibler (KL) divergence function will measure the disparity between the two distributions. Equation (1) describes the divergence function in terms of the VFE and a surprisal term - the negative log probability distribution of outcomes:

$$D_{KL}[q(x)||p(x|\phi)] = \mathbb{E}_{q(x)}\left[\ln\frac{q(x)}{p(x|\phi)}\right] =$$
$$\int q(x) \ln \frac{q(x)p(\phi)}{p(x,\phi)} dx = \int q(x) \ln \frac{q(x)}{p(x,\phi)} dx + \ln p(\phi) \quad (1)$$

where the VFE is defined as:

$$F \triangleq -\int q(x) \ln \frac{q(x)}{p(x,\phi)} dx \quad (2)$$

The above equations indicate that minimizing $F$ involves reducing the divergence between $q(x)$ and $p(x)$ while improving the model's predictive accuracy for $\phi$. Additional details are presented separately for a building agent and the community agent.

## III. Building Level AIF Agent Model

For a building continuous AIF model $g(.)$ includes a deterministic term and a noise term, $\varepsilon_z(t) \sim \mathcal{N}(0, \sigma_z^2)$, i.e.,

$$p(\phi(t)|x(t), a(t)) = g(x(t), a(t)) + \varepsilon_z(t) \quad (3)$$

where $x(t) = \begin{bmatrix} x_1(t) \\ \cdots \\ x_{|\chi|}(t) \end{bmatrix} = \mu_{prior}(t) + \varepsilon_\omega(t)$, $\mu_{prior}$ is prior data on state beliefs and $\varepsilon_\omega(t) \sim \mathcal{N}(0, \Sigma_\omega(t))$ and $\Sigma_\omega$ is the covariance function capturing uncertainty in the priors. The probability density of hidden states under this model is given by $p(x) = \mathcal{N}(\mu_{prior}, \Sigma_\omega)$. For presentation purposes, we will assume that a building AIF agent has two hidden states, namely, the occupancy level ($x_1 = OCC$) and the external infiltration airflow ($x_2 = \dot{m}_{ext}$). We note that air infiltration can be caused by an open window or door. Moreover, we will assume that $\phi(t)$, the indoor temperature, is the only sensory data from a building zone. For further simplification, we will assume that the example building has only one zone. These assumptions can be relaxed assuming that the building multiple zones function independently, and observations are statistically independent. The building control agent uses noisy measurement $\phi(t) \sim \mathcal{N}(T(t), \sigma_z^2)$ to refine its beliefs about the hidden states, i.e., occupancy level and infiltration airflow, where $\sigma_z^2$ the observation noise variance.

For illustration purposes and without loss of any generality, we will assume that the building AIF agent manages two control actions: the HVAC airflow ($a_1 = \dot{m}$) and the HVAC supply air temperature ($a_2 = T_{sup}$). These actions directly influence the zone's thermal dynamics and enable real-time adjustments to maintain occupant comfort and optimize energy usage. By modulating $\dot{m}$, the system can control the rate of air exchange while $T_{sup}$ determines the temperature of supplied air. These actions impact how quickly indoor temperatures reach desired setpoints. For example, during periods of high occupancy, increasing airflow or lowering supply air temperature can rapidly cool the space and improve occupant comfort. We now have

$$g(x(t), a(t)) = T(t+1) = T(t) + \Delta t \left[ \frac{\dot{m}(t)c_p(T_{sup}(t)-T(t))+OCC(t).q_{occ}+(U_w+\dot{m}_{ext}(t)c_p)(T_{amb}(t)-T(t))}{C_b} \right] \quad (4)$$

Here, $\dot{m}(t)$ is the HVAC airflow, $T_{sup}(t)$ is the supply air temperature, $OCC(t)$ is the occupant-driven internal heat generation, $\dot{m}_{ext}(t)$ is the infiltration airflow from outside, $c_p$ is the air's specific heat capacity, $U_w$ is the building envelope's heat transfer coefficient, $T_{amb}(t)$ is the outdoor temperature, and $C_b$ is the building's effective thermal capacity. Note that with expanded states and actions, $g(x(t), a(t))$ may involve additional terms and some of these terms may require approximations, but that discussion is outside the scope of this article. For the building agent $\mu_{prior} = [\mu_{OCC,prior}, \mu_{\dot{m}_{ext},prior}]$ encodes the expected occupancy profiles and infiltration patterns, which reflect the prior beliefs. If a dynamic model is available, such as a Markov process, the prior mean can be replaced with a state transition model:

$$x(t+1) = f_x(x(t)) + \varepsilon_\omega(t) \quad (5)$$

This dynamic function framework ensures that the system leverages historical data, time-of-day patterns, and other prior knowledge to infer hidden states effectively.

The building AIF agent also tries to shape future observations in accordance with some preferences. We introduce a target (or preferred observation) $\rho$ and define a preference distribution $p(\phi; \rho)$ that assigns higher probabilities to observations closer to $\rho$. This modifies the generative model to $p(\phi, x|a, \theta) = p(\phi|x, a, \theta) p(x|a, \theta)$. Although $\rho$ is a desired state rather than a random variable, it is treated as a prior preference in Active Inference. By incorporating this into the VFE equation and simplifying the KL divergence equation, the VFE equation becomes:

$$F = D_{KL}[q(x)||p(x)] - \mathbb{E}_{q(x)}[\ln p(\phi|x)] - \ln p(\phi; \rho) \quad (6)$$

Since $p(\phi; \rho)$ is independent of $x$, it simplifies to a constant term with respect to $q(x)$. Hence, the final VFE equation in the context of our building energy management problem becomes:

$$F = \frac{1}{2\sigma_z^2}[\phi(t) - T(t)]^2 + \frac{1}{2\sigma_\rho^2}[\rho(t) - \phi(t)]^2 + \frac{1}{2\sigma_\omega^2}\sum_{i \in \{OCC, \dot{m}_{ext}\}}[\mu_i(t) - \mu_{prior,i}(t)]^2 \quad (7)$$

Here, the first term describes the observation error, the second term encourages the indoor temperature to align with the target setpoint, and the last term penalizes deviations from the prior beliefs about hidden states. Minimizing this VFE function with respect to the hidden states $\mu_i(t)$ (occupancy level and infiltration) updates the aggregator's posterior states. This is formulated as $\mu_i^* = \arg\min_{\mu_i} F(\mu, a, \phi)$, and the iterative update for each hidden state is given by:

$$\mu_i(t) \leftarrow \mu_i(t) - \eta_i \frac{\partial F(t)}{\partial \mu_i(t)} \quad (8)$$

where $\eta_i$ is the learning rate.

On the other hand, minimizing VFE with respect to HVAC control actions ($a_1(t)$ and $a_2(t)$) optimizes the building's temperature control and guides the way toward the target setpoint $\rho(t)$ while minimizing energy costs. This is formulated as $a_i^* = \arg\min_{a_i} F(\mu, a, \phi)$, and the iterative update for controls is given by:

$$a_i^t \leftarrow a_i^t - \zeta_i \frac{\partial F}{\partial a_i^t} \quad (9)$$

where $\zeta_i$ is the learning rate.

The control goal in the building energy management problem is to adjust the HVAC airflow ($a_1$) and supply air temperature ($a_2$) to achieve the target indoor temperature ($\rho$) while minimizing energy use and VFE. By iteratively updating $a_1$ and $a_2$, the agent minimizes the free energy in order to align indoor temperature $\phi(t)$ with the target setpoint $\rho(t)$. This optimization ensures efficient energy consumption while maintaining occupant comfort. The building-layer AIF simultaneously infers hidden states and optimizes continuous HVAC actions. This combined approach allows the system to adapt to varying conditions and achieves robust performance in real-world building environments.

IV. COMMUNITY LEVEL DISCRETE AIF AGENT MODEL

For the illustration purposes, we will assume that the community manager oversees $N=2$ buildings, $K=2$ shared assets (photovoltaic (PV) panels and Energy Storage System (ESS)), manages the market Interactions, i.e., energy purchases or sales in the spot market, and carries out the Day-Ahead (DA) planning. The manager's primary objective is to minimize the community's net power costs while adhering to the DA energy purchase strategy. On the day of operation, the manager aims to minimize reliance on the more expensive spot market while accounting for penalties related to deviations from the DA plan, underutilization of PV generation, and inefficient battery discharge. Additionally, the manager strives to maintain an accurate knowledge of hidden states, such as building load states and the ESS State of Charge (SoC). To achieve these objectives, the aggregator employs a two-pronged strategy: (1) State inferencing using VFE (similar to what was presented above but for a discrete case), and (2) Policy selection using EFE, which combines a cost function with an entropy term to reduce state uncertainty. By penalizing ambiguity, the manager prioritizes actions that simultaneously optimize costs and improve state knowledge of buildings and common assets under its supervision, and external factors, such as energy market conditions. The community AIF agent implements a rolling horizon approach to balance short-term decision-making with long-term optimization. For instance, actions taken to reduce uncertainty about building states or battery SoC in one step improve the accuracy of subsequent decisions. This adaptability allows the framework to respond dynamically to real-time conditions, such as fluctuations in PV generation or occupant behavior.

The community AIF agent defines two types of hidden states in $\{s \in \mathcal{S}\}$, namely,

a. Building States, $(s_{B1}$ and $s_{B2}) \in$ {High Load, Med Load, Low Load}, and

b. ESS State of Charge, $s_{ESS} \in$ {E=Empty, L=Low, H=High, F=Full}.

PV states can be inferred from the states of ESS. These hidden states of buildings, when inferred, reflect their occupancy levels and infiltration airflow without violating the building's privacy. For instance, inferring a "High Load" state indicates to the manager that the building has high occupancy and/or potential energy waste due to open doors or windows. The combined hidden state space consists of 3×3×4=36 discrete states. The observation set $(o \in \mathcal{O})$ captures partial signals from buildings and the shared resource block:

a. Energy consumption level signals from buildings, $o_{B1}, o_{B2} \in$ {Up, Same, Down}, and

b. Electricity flow signals, $o_{ESS} \in$ {Import, Neutral, Export}.

By combining these signals, the model forms 3×3×3=27 observation categories. These partial observations update the manager's belief over the 36 hidden states and enable dynamic adjustments based on real-time information.

Actions, $\mathbf{u} \in \mathcal{U}$, are defined for buildings, ESS, and market transactions as follows:

a. Building actions are represented as $u_B \in$ {No Change, Small Reduction, Big Reduction},

b. Shared assets actions are $u_{ESS} \in$ {Charge, Hold, Discharge}, and

c. Market actions are denoted as $u_M \in$ {Buy, No Transaction, Sell}.

For simplicity, we assume the building actions are the same for both buildings. The combined action space comprises 3×3×3=27 discrete actions. The model structure allows the manager to balance energy usage, battery operation, and market transactions.

The discrete AIF model uses two key matrices to represent the relationship between observations, hidden states and transitions under different actions. We start with the observation matrix $A \in \mathbb{R}^{(num_{obs}) \times (num_{states})}$ which defines the probability of observing $o$ given hidden state $s$, i.e., $p(o|s) = A[o,s] \in \mathbb{R}^{27 \times 36}$. Here, columns represent hidden states and rows correspond to observations. $A$ is composed of the building observation submatrices ($A_{B1}$ and $A_{B2}$) and the ESS observation submatrix ($A_{ESS}$). The full observation matrix can be assembled by taking the Kronecker product of these submatrices:

$$A = A_{B1} \otimes A_{B2} \otimes A_{ESS} \qquad (10)$$

The transition matrices $B[\mathrm{u}(t)] \in \mathbb{R}^{36 \times 36}$ describe how hidden states evolve from one time step to the next under a given action $\mathrm{u}(t)$. Since the model includes 27 discrete aggregator actions, 27 distinct transition matrices exist. Each matrix corresponds to a specific action and captures the dynamics of both building states and ESS SoC. To construct these matrices, we must build the following submatrices: (a) Building transition submatrices ($B_{B1}$ and $B_{B2}$) which are 9×9 and explain how the building states evolve under different actions, and (b) ESS transition submatrix ($B_{ESS}$) which is 4×4 and explain the ESS SoC transitions and captures charging, discharging, or holding actions. The complete transition matrix for a given action is formed by combining these submatrices using the Kronecker product:

$$B = B_{B1} \otimes B_{B2} \otimes B_{ESS} \qquad (11)$$

Similar to continuous AIF, the VFE in discrete AIF measures how well an approximate variational posterior $q(s)$ aligns with the true posterior $p(s|o)$. In discrete AIF, this minimization results in an efficient belief update rule:

$$q(s') \propto q(s) B[u](s',s) A(o',s') \qquad (12)$$

This equation indicates that the aggregator updates its belief about the hidden state $(q(s'))$ based on prior beliefs $q(s)$, transition dynamics $B[u](s',s)$, and new observations $A(o',s')$. Maintaining an accurate posterior belief about hidden states (which are abstract and high level at the community level), enables more informed and adaptive decision-making. For instance, if building signals indicate a sudden increase in energy consumption, the manager can proactively request load reductions or discharge the battery to mitigate the surge to ensure responsive energy management.

For planning purposes, we will assume that the planning horizon is a full day, discretized into 15-minute intervals, resulting in ninety-six (96) time steps. The manager employs a rolling horizon strategy, such that at each step $t$, it solves for an optimal 4-step sub-horizon policy ($\pi_{t:t+3}$) over the next 45 minutes. This way, the model considers both the immediate and future impacts of its actions. This approach enables the system to adapt dynamically to new information, such as updated occupant signals or changes in electricity prices. At each time step the manager follows this sequence:

1. Horizon: Set the upcoming time steps ($t, \ldots, t + 3\Delta t$).
2. Sample Policies: Possible action sequences ($\pi = (u_t, u_{t+1}, u_{t+2}, u_{t+3})$) are generated.
3. Forward-Simulate Outcomes: For each policy, the agent predicts its impact on costs, energy usage, and uncertainty reduction using EFE.
4. Apply First Action: The first action of the selected policy is implemented, while the remaining horizon is re-evaluated at the next time step with updated observations.

This receding-horizon strategy ensures that decisions are made iteratively, with adjustments reflecting real-time conditions. Our formulation of EFE for a given policy includes epistemic and extrinsic terms. The epistemic term selects actions that improve state estimation and reduce uncertainty, while the extrinsic term ensures that selected actions align with prior preferences. We have:

$$G_\pi = D_{KL}\big(q(o^t|\pi^t)||p(o^t)\big) + \sum_s q(s^t|\pi^t) H\big(p(o^t|s^t)\big) \quad (13)$$

In Equation (13), the first term describes the risk and measures how the agent's predicted observations $q(o^t|\pi^t)$ deviate from preferred observations $p(o^t)$, and the second term shows the ambiguity which captures how uncertain future observations are, given each hidden state. Higher emphasis on ambiguity means the agent chooses future states that generate unambiguous and informative outcomes, while stronger preference demands can reduce exploration in favor of meeting cost or comfort objectives. For a 4-step horizon, the EFE for a policy $\pi_{t:t+3}$ can be rewritten with cost and ambiguity terms, as follows:

$$G_{\pi_{t:t+3\Delta t}} = \sum_t^{t+3\Delta t} \mathbb{E}_{q(s^t|\pi^t)}\Big[[Cost(o^t, u^t)] + \alpha_{amb} H\big(p(o^t|s^t)\big)\Big] \quad (14)$$

Equation (15) describes a cost function to align real-time community energy usage with a pre-defined day-ahead (DA) energy plan ($P_{target}^t$) while encouraging efficient use of on-site PV and battery storage and is defined as:

$$Cost(o^t, u^t) = Cost_{spot} + P_{dev}(|R^t|) + P_{UPV} + P_{Bat} \quad (15)$$

where $R^t = P_{RT}^t - P_{target}^t$.

By incorporating ambiguity penalties into the decision-making process, the manager prioritizes actions that reduce uncertainty about hidden states. This proactive approach improves decision quality over time as the system becomes more confident in its state estimates. Additionally, the integration of cost optimization—through spot market interactions, adherence to day-ahead plans, and efficient PV utilization—ensures the framework minimizes energy costs while balancing system constraints.

## V. SIMULATION RESULTS

In the simulations, each building is represented by a digital twin, namely an EnergyPlus model, which acts like the generative process and is capable of capturing heat transfer and airflow between inside and outside of a zone, occupants, building enclosure and windows, and energy assets. The simulated community consists of two buildings, each with an HVAC system, and a shared internal resource block with an ESS and PV panels. The digital twin integrates real-time factors such as energy usage, occupant schedules, PV generation, and battery dynamics. The community interacts with the energy market through day-ahead scheduling and real-time spot market adjustments. The simulations were done for reference residential buildings in Newark, NJ, for actual representative summer conditions. In the building-layer AIF, a given day is divided into 288 5-minute timesteps. Model parameters for building-layer agents are detailed in Table I, while community-layer parameters are summarized in Table II. Electricity prices were sourced from PJM datasets, and simulations were executed on an Apple MacBook Pro with an M4 processor, 16 GB of RAM, and a 10-core CPU.

TABLE I. BUILDING-LAYER AIF MODEL PARAMTERES

| Parameter | Value | Parameter | Value |
|---|---|---|---|
| Building thermal inertia | 2×10⁶ (J/°C) | Comfort band | ±2 (°C) |
| Specific heat capacity of air | 1005 (J/kg°C) | Generated heat by each occupant | 102 (W) |
| Min. and Max. HVAC supply air temperature | 10-25 (°C) | Min. and Max. HVAC airflow rate | 0 – 0.3 (kg/s) |

TABLE II. COMMUNITY-LAYER AIF MODEL PARAMTERES

| Parameter | Value | Parameter | Value |
|---|---|---|---|
| Maximum PV Panel Power Generation | 20 (kW) | Battery Cap. and Max. Power Output | 5 (kWh) - 20 (kW) |

First, we compare the AIF agents to a fully observable and deterministic optimization model which is widely regarded as the benchmark for decision-making in fully observable environments. We employed a MINLP approach, using the Gurobi solver with all equations and constraints specified in the building-layer AIF model. For the comparison, two types of AIF models were considered: a one-step look-ahead and a full-horizon AIF.

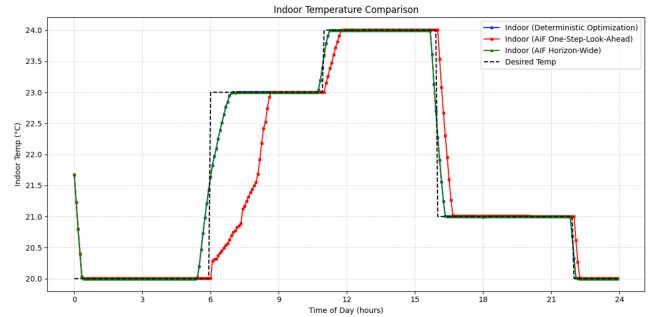

Figure 2. Comparison of predicted indoor temperatures for deterministic optimization and Active Inference models.

Figure 2 compares the predicted indoor temperature of the three models against the desired temperature setpoint. The deterministic optimization model and the horizon-wide AIF produce identical results, which demonstrate the ability of the horizon-wide AIF to replicate deterministic optimization when perfect information is available. In contrast, the one-step look-ahead AIF exhibits slight deviations during transitional periods, due to its shorter planning horizon. These deviations highlight an important trade-off in AIF models: while a horizon-wide AIF ensures long-term stability, a one-step look-ahead approach enhances short-term adaptability but may introduce minor fluctuations.

Figure 3 compares the AIF with the Deep Deterministic Policy Gradient (DDPG) algorithm, an actor-critic RL method suited for continuous action spaces like HVAC airflow and supply temperature control. DDPG learns optimal policies by interacting with the environment, observing states (indoor and ambient temperatures), and selecting corresponding actions. In contrast, AIF explicitly estimates hidden states, such as occupancy and infiltration, while optimizing temperature control. The RL model for this comparison consists of two hidden layers with 256 neurons with ReLU activation function. It was trained on one year of simulated data from EnergyPlus for the same reference building located in Newark, NJ. The dataset was split into 80% for training and 20% for testing. Both RL and AIF models aim to maintain room temperature within the desired range, but they differ in how they handle uncertainty and non-stationary environmental conditions. Unlike AIF, the RL model does not have direct access to hidden states, and it relies only on observed data and its prior training experience to react and predict indoor temperature. If, for example, occupant behavior or infiltration rates change abruptly, the RL controller may require retraining or additional adaptation mechanisms to maintain performance. AIF, on the other hand, can rapidly re-estimate hidden states and adjust accordingly.

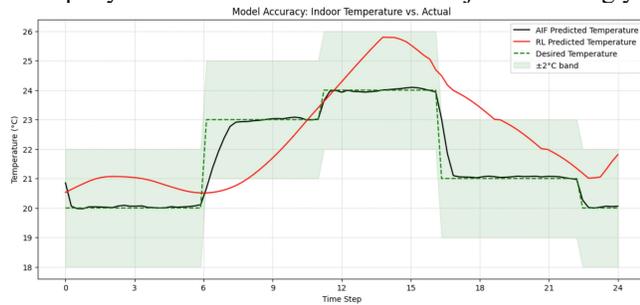

Figure 3. Comparison of the results of AIF with an RL algorithm for predicted indoor temperatures.

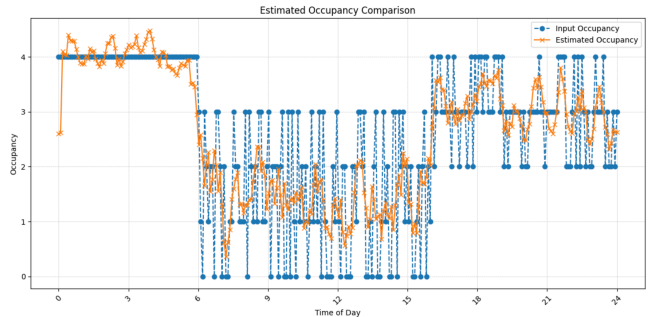

Figure 4. Comparison of the inferred occupancy profile against the original occupancy profile.

Figure 4 presents a comparison between the original occupancy profile used in the EnergyPlus simulation and the occupancy estimated by the AIF model over a 24-hour period. The AIF model effectively tracks fluctuations in occupancy, capturing both high and low activity periods with reasonable accuracy. Although some discrepancies exist, particularly during transitions between occupancy states, the model demonstrates an overall ability to adapt to dynamic occupancy patterns. This adaptive behavior enables the AIF model to make optimal real-time HVAC control decisions by continuously updating its internal beliefs in response to observed conditions.

At the community-layer, the AIF agent operates over a 24-hour horizon, which is divided into 15-minute intervals using a short-horizon EFE minimization. Figure 5 describes the inferred hidden states of the discrete AIF agent. For example, fluctuations between "High,High" and "Low,Med" illustrate the agent's changing beliefs about occupancy patterns, household appliance usage, HVAC power consumption, and/or other factors that drive building loads. It can also be seen that the states occasionally jump within short time spans. This indicates that the agent rapidly updates its hypotheses whenever it sees new signals.

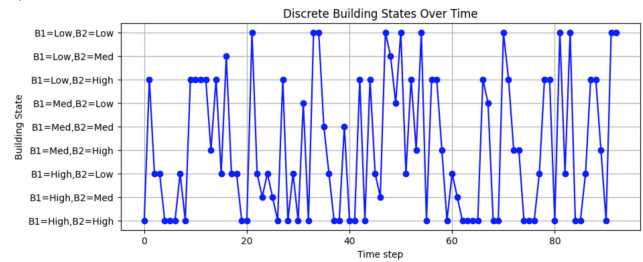

Figure 5. Inferred states of the buildings in the community.

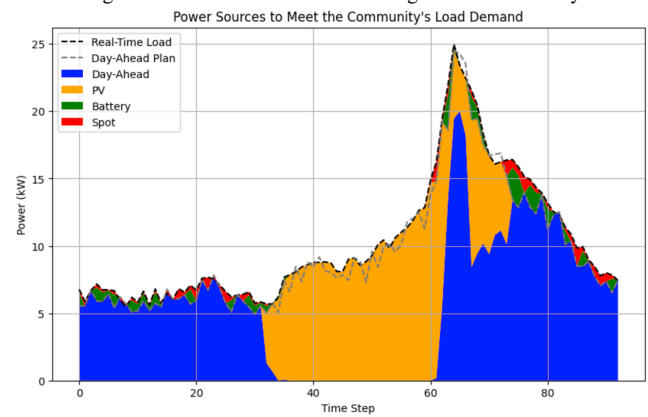

Figure 6. Real-time operation of the community energy management.

Figure 6 integrates real-time operations with day-ahead planning by visualizing the stacked power contributions from PV generation, battery adjustments, and spot market transactions, overlaid on the baseline day-ahead power purchase. Real-time operations closely follow the day-ahead plan, with minor deviations managed through dynamic adjustments in PV, battery, and market actions. The plot also highlights instances where the day-ahead forecasted demand exceeds real-time usage. In such cases, the system leverages excess day-ahead energy to either charge the battery or sell surplus power back to the spot market. When PV output is minimal at the beginning of the operation day, the system predominantly relies on the spot market and battery storage to meet residual demands. As PV generation ramps up throughout the day, it becomes the primary energy source and reduces the dependency on external sources. During peak demand, the combined utilization of PV generation and battery storage effectively reduces reliance on the spot market. Later in the day, as demand decreases and PV generation subsides, the system shifts its focus to discharging the battery to reduce spot market transactions. This visualization highlights the framework's ability to balance multiple energy sources, adapt to real-time uncertainties, and optimize resource utilization for improved cost efficiency.

Figure 7 illustrates the impact of incorporating ambiguity into the stepwise EFE values for five values of the ambiguity coefficient, $\alpha_{amb}$ (0.0 to 2.0). Increasing $\alpha_{amb}$ places greater emphasis on reducing ambiguity, which generally results in higher stepwise costs. However, models with lower ambiguity can occasionally experience higher costs at specific timesteps, depending on how the environment responds to uncertainty. Interestingly, agents with higher $\alpha_{amb}$ can sometimes select more effective policies that resolve state uncertainty while also reducing overall costs, resulting in lower cumulative EFE. These exploratory actions incur higher short-term expenses but enhance state estimation and long-term decision-making. However, when environmental feedback offers limited benefits within the simulation horizon, the added costs of minimizing ambiguity may not be recovered. Thus, while ambiguity minimization can improve outcomes, it does not universally guarantee lower total costs.

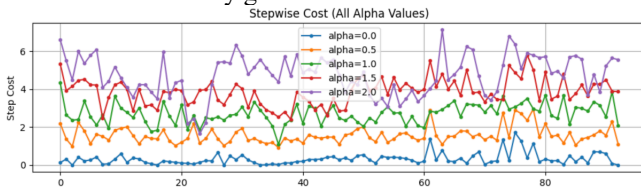

Figure 7. Impact of ambiguity term on the Expected Free Energy.

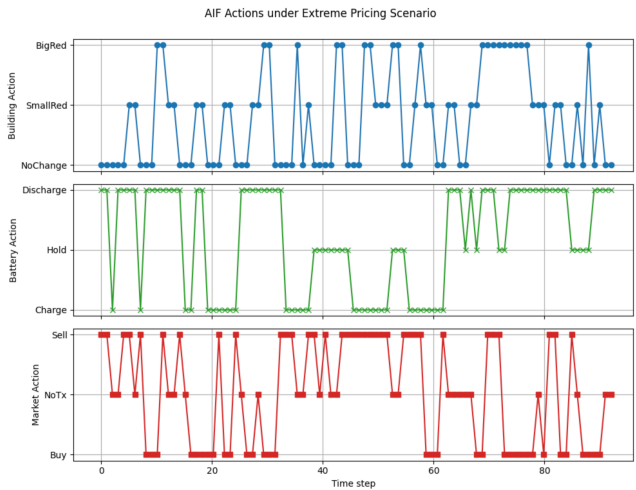

Figure 8. The aggregator actions under the extreme pricing scenario.

Finally, we evaluate the AIF performance under an extreme real-time pricing (RTP) scenario. We consider a summer day with a power generation shortage and elevated electricity costs. In this scenario, RTP rates were doubled, with peak-hour prices (16:00–20:00) reaching three times the baseline, resulting in higher utility bills for end-users. To mitigate these costs, the community AIF controller optimizes energy management by reducing reliance on the high-priced spot market. Figure 8 illustrates the system's decision-making for three action types. The model monitors real-time energy consumption in buildings and issues demand reduction requests based on grid conditions and pricing signals. It also discharges the battery during peak hours to reduce grid dependence and preserves or recharges it during lower-cost periods. Moreover, the model determines market transactions and decides when to buy, sell, or abstain from trading based on real-time conditions. This strategy led to decreased power purchases during peak pricing and increased energy sales from battery discharge and surplus PV generation. Overall, the results highlight the model's effectiveness in optimizing battery usage and market transactions.

## VI. CONCLUSION

This paper introduced a dual-layer Active Inference architecture for building community energy management, marking the first application of AIF principles to planning and control in energy systems. Two distinct AIF models were developed: a continuous framework to control HVAC operations at the building level and a discrete framework to coordinate power flows, battery usage, and external market interactions at the community level. Both layers try to infer their respective hidden states from observations and ensure that the performance targets are tracked and achieved. At the community level, the ambiguity term plays a significant role in balancing uncertainties with achieving short-term benefits. The results from the proposed AIF models compared to exact optimization under full observations and with RL under similar partial observations. The results show the same or superior performance by the AIF framework. The formulation is quite generic and, under some mild assumptions, can be extended to communities with a larger number of buildings and zones within buildings.